\documentclass[11pt]{article}
\usepackage{a4}
\usepackage{amssymb}
\usepackage{hyperref}

\let\nonu=\nonumber

\def\su{\ifmmode SU(2) \else $SU(2)$\fi}
\def\sb{{\bar s}}
\def\tb{{\bar t}}

\begin{document}

\title{On the icosahedron: from two to three dimensions.}
\author{Marc P. Bellon\thanks{On leave from:
Laboratoire de Physique Th\'eorique et Hautes Energies, Boite 126,
4 Place Jussieu,  75252 Paris Cedex 05. Unit\'e Mixte de Recherche
UMR 7589, Universit\'e Pierre et Marie Curie-Paris6; CNRS;
Universit\'e Denis Diderot-Paris7.}\\
 {\normalsize \it CEFIMAS, Av.\,Santa Fe 1145,
 1069 Capital Federal, Argentina}\\
\normalsize\it Departamento de F\'\i sica,
Universidad Nacional de La Plata\\ {\normalsize\it C.C. 67, 1900
La Plata, Argentina}}
\date{}
\maketitle
\begin{abstract}
 In his famous book~\cite{klein}, Felix Klein describes a complex
variable for the quotients of the ordinary sphere by the finite
groups of rotations and in particular for the most complex situation
of the quotient by the symmetry group of the icosahedron. The
purpose of this work and its sequels is to obtain similar results
for the quotients of the three--dimensional sphere. Various
properties of the group $SU(2)$ and of its representations are used
to obtain explicit expressions for coordinates and the relations
they satisfy.
\end{abstract}

\section{Introduction}
The study of the quotients of the two dimensional sphere by subgroups of the
rotation group has a long story. The identification of the sphere with the complex
projective sphere has linked this study to the one of binary forms and their
invariant combinations (see in particular~\cite{WBGF}), and it has been masterfully
exposed in the work of Felix Klein~\cite{klein} more than a
century ago. It has nevertheless seen a renewed interest with its
application to the classification of singularities of surfaces
through the McKay correspondence~\cite{McK2,McK3,GV}.

The quotients of the three--dimensional sphere are characterized by the
same subgroups of \su\ and were classified long ago~\cite{TS30}. The
possible cosmological interest of such quotients, as proposed by Luminet
et al.~\cite{luminet}, renewed the interest for the functions defined on
such quotients. In recent works, a description of the functions on the regular quotients
of the three sphere has been reached~\cite{We05,Be06}, but the algebra formed by
these functions was not considered.  This work aims at filling this gap, describing
generators of these algebras and the relations they satisfy. The reverse problem of
identifying from a description by generator and relations the algebra of functions
of some quotient of the three sphere will be reserved for a subsequent
work~\cite{Be07}.

This work uses the fundamental fact that the functions on the quotients
form multiplets of \su. When considering
polynomials in such functions, the use of Clebsch--Gordan coefficients
allows to group them in irreducible representations of \su. This produces
either new multiplets of functions or relations.

The next section will recall the procedure for the enumeration of the functions on
the quotients. Then the main technical tool is introduced, which allows to relate
operations in the three-dimensional realm to known invariant combinations in the
two-dimensional one. This allows to give a general description of the generators of
the algebra of functions. In some cases, one can use sections of non-trivial line
bundles to obtain embeddings of the quotient in a projective space. Finally,
I show how this structure is realized for every of the finite subgroups of \su.

\section{The functions.}
The functions on the quotients of the
three--sphere, classified by eigenspaces of the Laplacian, can be determined
from the special structure of the three--sphere. Indeed, this sphere can be given
the group structure of \su\ or the unit quaternion. Its symmetry group is the product
of its right and left actions on itself. Functions on the
sphere form a representation of the product of these two \su\ groups. Its
decomposition in irreducible representations involves all the possible
product of a representation of $\su^L$ with the equivalent one of $\su^R$.

The
identification with \su\ is made clear by introducing complex variables:
\begin{equation}
    \label{su2}  p = \pmatrix{s&-\bar t\cr t& \bar s\cr}
\end{equation}
with the constraint $s \bar s + t \bar t=1$.  With these variables, the
action of the Lie generators of the groups $\su^L$ and $\su^R$ can be
written as:
\begin{eqnarray}
     J_z^L &=& {\textstyle {1\over2}} (s\partial_s + \tb\partial_\tb
        - t \partial_t -\sb \partial_\sb), \nonu\\
     J_+^L &=& s\partial_t - \tb \partial_\sb, \nonu\\
     J_-^L &=& t\partial_s - \sb \partial_\tb, \label{jL}\\
     J_z^R &=& {\textstyle {1\over2}} (s\partial_s - \tb\partial_\tb
        + t \partial_t -\sb \partial_\sb), \nonu\\
     J_+^R &=& -s\partial_\tb + t \partial_\sb, \nonu\\
     J_-^R &=& -\tb\partial_s + \sb \partial_t.\label{jR}
 \end{eqnarray}
The quotients we consider in this work are quotient by a subgroup of one
of the two \su\ subgroups, so that they are automatically regulars. If we
suppose that the monodromy group is a subgroup of $\su^L$, the
functions on the quotients will appear in multiplets of $\su^R$ and new
functions can be obtained from a given one through the action of $J_-^R$
or $J_+^R$. The highest weight vector in a
representation of $\su^R$ can be represented by a function of only $s$ and $t$.

The functions on the quotients of the three-dimensional sphere by a
subgroup $\Gamma$ of $\su^L$ can therefore be described as follows. They form
multiplet of $\su^R$ in which all elements can be deduced from the highest weight
vector by the action of $J_-^R$. The possible highest weight vectors are
functions of $s$ and $t$ invariant under $\Gamma$. The invariant functions are
known since the work of
Klein~\cite{klein}. They are generated by three functions satisfying an
algebraic relation and can all be made from the lowest degree one by the
action of the Hessian and the cross-product. In the sequel, we show how
to go from this generation through the action of the lowering operator
$J_-^R$ to a direct generation by algebraic operations.

\section{Clebsch--Gordan }
Since the functions come in multiplets of \su, the different products of
functions in multiplets of spin $j_1$ and $j_2$ of \su\ belong to the
tensor product $(j_1) \otimes (j_2)$, which can be decomposed in the sum
of representations $(j_1+j_2)$ down to $(|j_1-j_2|)$.

We model the representations on the one by homogeneous polynomials of
degree $2j$ in two variables $u$ and $v$. The monomials $\langle j,m\rangle =
u^{j+m} v^{j-m}$
are vectors with the eigenvalue $m$ of $L_z$, but they are not normalized.
However the representations of $L_+$ and $L_-$ in this basis is simple
with integer matrix elements:
\begin{equation}
    L_+( \langle j,m\rangle ) = (j-m)\langle j,m+1\rangle,
      \qquad L_-( \langle j,m\rangle) = (j+m)\langle j,m-1\rangle.
\end{equation}
This translates in simple formulas for
the Clebsch--Gordan coefficients. Differential formulas have been given
in~\cite{Be06} for the general case and do not involve any squared
roots. In this work, we are mainly interested in the highest
weight vectors of the representations appearing in the composition of two
irreducible representations. We have, up to normalization:
\begin{equation} \label{comb}
    \langle j_1+j_2-k,j_1+j_2-k\rangle = \sum_{l=0}^k (-1)^l \pmatrix{l\cr k}
    \langle j_1,j_1-l\rangle \langle j_2,j_2-(k-l)\rangle,
\end{equation}
since $L_+$ annihilates the right hand side of this equation.

In the case of the functions on the sphere, the element $\langle j,j\rangle$ of a
$\su^R$ multiplet is a homogeneous polynomial of degree $2j$ in $(s,t)$
and the other elements of the multiplet are obtained by the action of $L_-^R$.
\begin{equation} \label{lowering}
    \bigl(L_-^R\bigr)^l \langle j,j\rangle = {(2j)!\over (2j-l)!} \langle j,j-l\rangle
\end{equation}
In equation~(\ref{comb}), we need the vectors $\langle j_1,j_1\rangle$ up to
$\langle j_1,j_1-k\rangle$
and we therefore will multiply this equation by a factor
$(2j_1-l)!/(2j_1-k)!$ to have the common factor $(2j_1)!/(2j_1-k)!$. This additional
factor can be obtained by a degree counting operator, since the left hand side of
equation~(\ref{lowering}) is homogeneous of degree $2j-l$ in the
variables $s$ and $t$. The homogeneous degree is obtained with the
differential operator $s\partial_s+t\partial_t$. Switching to alternative variables
$s_1$ and $t_1$ and then back, the full factor is
obtained through powers of the operator $D_1=s\partial_{s_1}+t\partial_{t_1}$.
%applied to the function with $s$ and $t$ substituted by $s_1$ and $t_1$.
% it is advisable to express the function $\langle j_1,j_1\rangle$ in theand the function $\langle j_2,j_2\rangle$ in the variables $s_2$ and $t_2$. and in the other by $L_2 = \sb \partial_{t_2} - \tb \partial_{s_2}$. Finally, we write:
Finally, the change to the variables $s_1$ and $t_1$  can be done before the actions
of $L_-^R$ if we substitute this operator by $L_1=\sb \partial_{t_1} - \tb
\partial_{s_1}$. We thus obtain the following expression for the vectors appearing in
eq.~(\ref{comb}):
\begin{equation}
{(2j_1)! \over (2j_1-k)!} \langle j_1,j_1-l\rangle = (D_1^{l-k} L_1^k \langle j_1,j_1\rangle
)\bigr|_{s_1=s,t_1=t}
\end{equation}
Using the corresponding representation for the elements $\langle
j_2,j_2-(l-k)\rangle$, with the operator $L_2$ and $D_2$ defined in terms of
$\partial_{s_2}$ and $\partial_{t_2}$, we transform formula~(\ref{comb}) into:
\begin{eqnarray}    \label{comb2}
   {(2j_1)! (2j_2)! \over (2j_1-k)! (2j_2-k)!}\mskip -20mu
   &&\langle j_1+j_2-k, j_1+j_2-k \rangle
   \nonu\\
   &&= \sum_{l=0}^k (-1)^l
   \pmatrix{l\cr k} D_1^{l-k} L_1^k \langle j_1,j_1 \rangle
                D_2^k L_2^{l-k} \langle j_2,j_2 \rangle.
\end{eqnarray}
In this formula, $\langle j_i,j_i \rangle$ is understood as a function of $s_i$ and
$t_i$ and we left implicit the step of changing all remaining $s_i$ and $t_i$ variables
to $s$ and $t$. Now the differential operators $L_2$ and $D_2$ do
not act on the variables in $\langle j_1,j_1 \rangle$ and their
coefficients do not depend on the differentiation variable. This is also
true of $L_1$ and $D_1$ so we are
free to reorder the terms in~(\ref{comb2}). The sum on $l$ is therefore
the binomial expansion of $(D_1 L_2 - L_1 D_2)^k$ applied
to a product of polynomials in $(s_1,t_1)$ and $(s_2,t_2)$. Now
\begin{equation}\label{fond}
D_1 L_2 -L_1 D_2 = (s\sb + t \tb) (\partial_{s_1} \partial_{t_2} -
\partial_{t_1} \partial_{s_2}).
\end{equation}

We obtain therefore that in the combination of multiplets of spin $j_1$
and $j_2$ with the spin $j_1+j_2-k$ appears a factor $(s\sb+ t\tb)^k$,
which is 1 on the sphere. The same result could be obtained from the
condition that the highest weight in the multiplet of spin $j_1+j_2-k$
is annihilated by $L_+^R$, but we obtained a more precise result.
Expanding the result of equation~(\ref{fond}) and returning to the
only variables $s$ and $t$ shows that the highest weight vector can
be obtained by
\begin{equation}
 	{(2j_1-k)! (2j_2-k)!\over (2j_1)! (2j_2)! }\sum_{l=0}^k (-1)^l \pmatrix{k\cr l}
    (\partial_s^{k-l} \partial_t^l \langle j_1,j_1 \rangle )
    (\partial_s^l \partial_t^{k-l} \langle j_2,j_2 \rangle).
\end{equation}
This is the precise definition of the covariants $(f,g)^k$ used in the studies of
the nineteenth century~\cite{WBGF}, even the global factor is the same. This
formula was used for the combination of representations of \su\ defined as binary
homogeneous polynomials in~\cite{Be06}.

\section{Generators of the algebra of functions.}
Now the description of the algebra of functions on the quotient becomes
clear. We take a first multiplet $M$ of functions in the representation
of spin $j$ of $\su^R$, with a highest weight vector $P$. Their binary
products belong to the representations of spin $2j$, $2j-2$,\dots
The multiplet with spin $2j-2$ has $Q$, the Hessian of $P$ as highest weight,
from the result of the preceding section.  Among the ternary products of the
functions of $M$ we can select the representation of spin $3j-3$ of
$\su^R$, which can be obtained as a special component of the composition
of the multiplet of spin $j$ and $2j-2$. It will have as highest weight
$R$, the cross-product of $P$ and $Q$.

Klein proved that every invariant polynomial can be
expressed as a product of the polynomials $P$, $Q$ and $R$. Every multiplet
can therefore be obtained as the highest spin component in
the combination of $p$ times the multiplet $M$, $q$ times the multiplet
of spin $2j-2$ generated by $Q$ and $r$ times the multiplet generated by
$R$. Since the multiplets generated by $Q$ and $R$ are formed from
polynomials of degree 2 and 3 of the elements of the multiplet $M$, every
multiplet appears as a polynomial in the elements of $M$, which are therefore
coordinates for the quotient space. In the cases of the group $C_n$ and $D_{4n}$, this description is
not sufficient, since the first invariant polynomial $P=st$ has Hessian 1
in the case $C_n$ and
in the case $D_{4n}$, $P=s^2t^2$ has a Hessian proportional to itself.
The multiplet generated by $Q$, with a degree proportional to $n$ has to
be introduced independently.

Is is however not always necessary to use multiplets generated by invariants.

Indeed, let us consider a polynomial $V$ in $s$ and $t$ which belongs to some
one-dimensional representation of $\Gamma$. Since the actions of $\su^R$ and
$\Gamma$ commute, the polynomials generated from $V$ by the action of $L_-^R$
are multiplied by the same factor under the action of $\Gamma$. The multiplet they
form is not invariant under $\Gamma$, but the different transformations under
$\Gamma$ are projectively equivalent. We obtain an
imbedding in a projective space of smaller dimension than the affine space
whose coordinates are the invariant functions.

The definition of the quotient as an algebraic variety does not end with the
identification of coordinates. We further need to describe the relations
they satisfy. Every combination of the coordinates with a
$\su^R$ spin which do not correspond to any invariant of $\Gamma$ must
vanish. This provides for a rich set of relations. The difficult part will be
to identify a sufficient set of relations from which all others can be
deduced. This will be done in~\cite{Be07}.

\section{Explicit results.}

\subsection{Cyclic group $\mathbb Z_n$.}
In this case, all representations of $\Gamma$ are one-dimensional, so that
we could use the doublets $(s,-\tb)$ and $(t,\sb)$ as projective coordinates.
However, this does not give an embedding of the quotient, since the map is
the Hopf fibration of the 3-sphere with image the 2-sphere represented as the
complex projective line. If we start from the lowest degree invariant, we
have three coordinates in the spin 1 representation of \su, which satisfy a
quadratic relation and correspond to the euclidian representation of the
two-dimensional sphere, again the quotient of the three-dimensional sphere by the Hopf fibration.

The $n$ dependent part comes from the second independent invariant $s^n+t^n$.
It gives rise to a representation of spin $n/2$. The combination with the coordinates of the two-sphere in a spin~1 representation yield representations of spin $n/2-1$, $n/2$ and~$n/2+1$. The representation of spin $n/2-1$ must be zero and correspond to constraints and the combination of spin $n/2$ gives the multiplet
generated by $s^n-t^n$. In the case where $n$ is odd,
the representations of half integer spin of \su\ are pseudoreal and
complex conjugation relates the multiplet generated by $s^n+t^n$ to the one
generated by $s^n-t^n$. This means that the reality conditions relate the
multiplet of spin $n/2$ and its combination with the 3 coordinates of the
two-sphere.

\subsection{Dihedral groups $D_{4n}$.}

This case is completely similar to the preceding one, apart from the fact
that $st$ and $s^n+t^n$ ($s^n+i t^n$ in the case where $n$ is odd) are in
non-trivial representations of $D_{4n}$. We
therefore have a description as a part of the product of the 2-dimensional
projective space and a $n$-dimensional projective space.  The relations
remain the same.

\subsection{Tetrahedral group.}

In this case, we have a polynomial of degree 4 which belongs to a non-trivial
one-dimensional representation of the group.
Klein introduced it as the
polynomial with zeros on the vertices of the tetrahedron:
\begin{equation}
    V = s^4 + t^4 + 2i\sqrt3 \,s^2t^2. \label{tetra}
\end{equation}
It is not fully invariant, being multiplied by a cubic root of the unity
under the action of $\Gamma$.
It is part of a representation of spin 2 of $\su^R$, which is
of dimension 5, providing an embedding of the quotient space in 4 dimensional
projective space. The lowest degree invariant polynomial is of degree 6:
\begin{equation}
    P= st(s^4-t^4), \label{tetra2}
\end{equation}
It will generate a representation of spin 3 of dimension 7.

The lowest dimensional embedding translates in a much simpler description.
Indeed the space of homogeneous polynomials of degree
$n$ has dimension proportional to $n^4$ in one case and to $n^6$ in the
other. When decomposing this space in irreducible representations of $\su^R$,
we have on the order $n$ representations of dimensions of order $n$ so that
the typical multiplicity of a representation is proportional to $n^2$
in one case and $n^4$ in the other.  Identifying a fundamental set of
representations becomes much more easy in the projective case.

Quadratic polynomials in a spin 2 representation decompose in spin 4,
spin 2 and spin 0 representations. The spin 0 representation can be seen to
be zero, while the spin 2 representation has a generator similar to the one
of $V$ apart from the substitution of $i$ by $-i$. The spin 0 representation
gives one quadratic relation on the coordinates. It therefore seems that we
simply recovered the three dimensional sphere as a quadric in a four
dimensional projective space, but in fact, the non trivial geometry comes
from the reality conditions, which link the fundamental spin 2 with the spin
2 part in the quadratic polynomials. This is coherent with the fact that the
projective coordinates are multiplied by a cubic root of unity under the action of
the group $\Gamma$.

The cubic polynomials in the projective coordinates are fully invariant
under $\Gamma$ and decompose into the representation $j=0$, $j=2$,
$j=3$, $j=4$ and $j=6$. The $j=2$ representation vanishes since it can be
written as the spin 0 relation times the fundamental representation. The
$j=0$ polynomial gives a number which can be used to normalize the
projective representation. The ambiguity in the projective coordinates
are reduced from the multiplication by any non zero complex number to
the multiplication by a cubic root of unity. The $j=3$, $j=4$ and $j=6$
representations are descendents of the three fundamental invariants,
which are of degree 6, 8 and~12.

Since all invariants of the tetrahedron group can be obtained as products
of the three fundamental ones, every function on the quotient of the
three--sphere can be obtained from a polynomial of degree $3n$ in the five
projective coordinates.

\subsection{Octahedron group.}
In this case also, we can also work with a projective representation. The
invariant polynomial $P$ of the tetrahedral group defined in~(\ref{tetra2}) is invariant
up to a sign under the symmetry group of the octahedron. It generates a
spin~3 representation of $\su^R$ which will give embedding in a six
dimensional projective space. In this case, the reality conditions can be
satisfied inside the multiplet, so that we are in a real projective space.

The quadratic polynomials decompose into components of spin~6, 4, 2 and~0.
The spin~6 and~4 components are the descendents of the first two invariants
of the group, the spin~2 one is the set of relations and the spin~0 objects
allows for the normalization of the projective coordinates.

Among the cubic polynomials appears a spin~6 representation generated by the third invariant of the tetrahedron group. The descendants of the third invariant of the
octahedron group form a representation of spin~9 which appears first the polynomials of degree four in the projective coordinates.

\subsection{Dodecahedron group.}

In this case, we cannot find a non trivial one-dimensional representation. We
therefore must start from the 13 fonctions in the multiplet of spin 6 generated
by the fundamental invariant of degree twelve of the dodecahedron group.

As in the other cases, among the quadratic polynomials we can find the
spin~10 representation generated by the invariant of degree twenty and the cubic
polynomials have a part of spin 15 which is generated by the invariant of degree
thirty. We therefore identified a set of coordinates, such that every function on
the Poincar\'e sphere is a function of these coordinates. The difficulty lies in
the identification of a fundamental set of equations satisfied by these
coordinates. The quadratic functions belong to the representations
of all even spin from 0 to 12. The absence of invariants of the dodecahedron
group means that the spin 2, 4 and~8 representations must be zero. The spin 6
representation in the quadratic part must be related to the fundamental spin 6
representation. In fact, we can use only part of these quadratic relations to generate all of them. This will be detailed in our next work~\cite{Be07}.

\section{Conclusion.}

I have shown that the algebra of functions on the quotients of the sphere
is generated by one or two \su\ multiplets of functions. The essential
ingredient is the use of Clebsch--Gordan coefficients to combine products in
multiplets: the relation~(\ref{comb}) links such operations on the functions
of the three-sphere to equivalent ones on the holomorphic sections of bundles
on the two--sphere, reducing to the situation studied by Klein.

If it is easy to identify relations from the multiplets which do not
correspond to any invariant of $\Gamma$, it is more difficult to identify a
sufficient sets of relations. In a companion paper, I will show that the
quadratic relations are sufficient to generate all relations.

\noindent{\bf Acknowledgement}: the French CNRS is the only source of financial support for this work through a ``mise \`a disposition". I am grateful to the Instituto Argentino de Matem\'atica (IAM-CONICET) and its director Pr.\ Gustavo Corach for its hospitality.

\end{document}